\documentclass[12pt]{amsart}
\usepackage{epsf, amssymb, graphicx, multicol, amscd,times,stmaryrd, amsmath}


\newtheorem{theorem}{Theorem}[section]

\newtheorem{fact}[theorem]{Fact}

\newtheorem{proposition}[theorem]{Proposition}
\newtheorem{conjecture}[theorem]{Conjecture}

\theoremstyle{remark}
\newtheorem{remark}[theorem]{Remark}
\newtheorem{example}[theorem]{Example}
\theoremstyle{example}
\theoremstyle{definition}
\newtheorem{definition}[theorem]{Definition}

\newcommand{\R}{\mathbb{R}}
\newcommand{\Z}{\mathbb{Z}}

\newcommand{\ve}{\varepsilon}

\numberwithin{equation}{section}

\title{A note on Lagrangian cobordisms between Legendrian submanifolds of $\R^{2n+1}$}

\author{Roman Golovko}

\address{D\'{e}partement de Math\'{e}matiques\\ Universit\'{e} Libre de Bruxelles\\\newline
CP 218, Boulevard du Triomphe, B-1050 Bruxelles, Belgium}
\email{r.a.golovko@gmail.com}

\date{\today}

\begin{document}
\begin{abstract}
We study the relation of an embedded Lagrangian cobordism between two closed, orientable
Legendrian submanifolds of $\R^{2n+1}$. More precisely, we investigate the behavior of the Thurston-Bennequin number and (linearized) Legendrian contact homology under this relation. The result about the Thurston-Bennequin number can be considered as a generalization of the result of Chantraine which holds when $n=1$. In addition, we provide a few constructions of Lagrangian cobordisms and prove that there are infinitely many pairs of exact Lagrangian cobordant and not pairwise
Legendrian isotopic Legendrian $n$-tori in $\R^{2n+1}$.
\end{abstract}
\maketitle

\section{Introduction}
\subsection{Basic definitions}
A \emph{contact manifold} $(M,\xi)$ is a $(2n+1)$-dimensional manifold $M$ equipped with
a smooth maximally nonintegrable hyperplane field $\xi\subset TM$, i.e., locally $\xi = \ker \alpha$, where $\alpha$ is
a $1$-form which satisfies $\alpha\wedge (d\alpha)^{n}\neq 0$. $\xi$ is a \emph{contact structure} and $\alpha$ is a \emph{contact $1$-form} which
locally defines $\xi$. The \emph{Reeb vector field $R_{\alpha}$} of a contact form $\alpha$ is uniquely defined by the conditions
$\alpha(R_{\alpha}) = 1$ and $d\alpha(R_{\alpha}, \cdot)= 0$.  The most basic contact manifold is $(\R^{2n+1}, \xi)$, where $\R^{2n+1}$ has coordinates $(x_{1}, y_{1},\dots, x_{n}, y_{n}, z)$, and $\xi$ is given by $\alpha=dz-\sum_{i=1}^{n}y_{i}dx_{i}$. Note that $R_{\alpha}=\partial_z$. From now on, for ease of notation, we write $\R^{2n+1}$ instead of $(\R^{2n+1},\xi)$.

A \emph{Legendrian submanifold} of $\R^{2n+1}$ is an $n$-dimensional submanifold $\Lambda$
which is everywhere tangent to $\xi$, i.e., $T_{x} \Lambda \subset \xi_{x}$ for every $x\in \Lambda$.
The \emph{Lagrangian projection} is a map $\Pi:\R^{2n+1}\to\R^{2n}$ defined by
\begin{align*}
\Pi(x_{1}, y_{1},\dots , x_{n}, y_{n}, z)=(x_{1}, y_{1},\dots, x_{n}, y_{n}).
\end{align*}
Moreover, for $\Lambda$ in an open dense subset of all Legendrian submanifolds with $C^{\infty}$ topology,
the self-intersection of $\Pi(\Lambda)$ consists of a finite number of transverse double points.
Legendrian submanifolds which satisfy this property are called \emph{chord generic}.
A \emph{Reeb chord} of $\Lambda$ is a path along the flow of the Reeb vector
field which begins and ends on $\Lambda$. Since $R_{\alpha}=\partial_z$, there is a one-to-one correspondence between Reeb
chords of $\Lambda$ and double points of $\Pi(\Lambda)$. From now on we assume that all Legendrian submanifolds of $\R^{2n+1}$ are connected and chord generic.

The \emph{symplectization} of $\R^{2n+1}$ is the symplectic manifold
$(\R\times\R^{2n+1}, d(e^{t}\alpha))$, where $t$ is a coordinate on $\R$.

\begin{definition}\label{scdefinition}
Let $\Lambda_{-}$ and $\Lambda_{+}$ be two Legendrian submanifolds of $\R^{2n+1}$.
We say that that $\Lambda_{-}$ is
cobordant to $\Lambda_{+}$ if
there exists a smooth cobordism $(L; \Lambda_{-},
\Lambda_{+})$, and an
embedding from $L$
to $(\R\times \R^{2n+1}, d(e^t\alpha))$ such that
\begin{align*}
L|_{(-\infty, -T_{L}]\times \R^{2n+1}}& = (-\infty, -T_{L}]\times \Lambda_{-},\\
L|_{[T_{L},\infty)\times \R^{2n+1}}& = [T_{L},\infty)\times \Lambda_{+}
\end{align*}
for some $T_{L}\gg 0$ and $L^{c}:=L|_{[-T_{L}-1,T_{L}+1]\times \R^{2n+1}}$ is compact. In the case of a Lagrangian (exact Lagrangian) embedding, we say that $\Lambda_{-}$ is Lagrangian (exact Lagrangian) cobordant to $\Lambda_{+}$.
We will in general not distinguish between $L$ and $L^{c}$ and call both $L$.
\end{definition}
 We now define the following notations. If $L$ is an embedded, an embedded Lagrangian or an embedded exact Lagrangian cobordism from $\Lambda_{-}$ to $\Lambda_{+}$, then we write $\Lambda_{-}\prec_{L} \Lambda_{+}$, $\Lambda_{-}\prec^{lag}_{L}\Lambda_{+}$ or $\Lambda_{-}\prec^{ex}_{L}\Lambda_{+}$, respectively. If $L_{\Lambda}$ is a filling, a Lagrangian filling or an exact Lagrangian filling of $\Lambda$ in the symplectization of $\R^{2n+1}$, i.e., $L_{\Lambda}$ is an embedded, an embedded Lagrangian or an embedded exact Lagrangian cobordism with empty $-\infty$-boundary and $+\infty$-boundary $\Lambda$, then we write $\emptyset\prec_{L_{\Lambda}}\Lambda$, $\emptyset\prec^{lag}_{L_{\Lambda}}\Lambda$ or $\emptyset\prec^{ex}_{L_{\Lambda}}\Lambda$, respectively.

For the discussion about Lagrangian cobordisms between Legendrian knots we refer to~\cite{Chantraine} and~\cite{EkholmHondaKalman}, and for the
obstructions to the existence of Lagrangian cobordisms defined using the theory of generating families we refer to~\cite{SabloffTraynor,SabloffTraynor1}.

\subsection{Legendrian contact homology}

Legendrian contact homology was introduced by
Eliashberg, Givental and Hofer in~\cite{EliashbergGiventalHofer} and independently, for Legendrian knots in
$\R^3$, by Chekanov~\cite{Chekanov}.
We now briefly remind the reader of the definition of linearized Legendrian contact homology complex of a closed, orientable, chord generic Legendrian submanifold $\Lambda\subset \R^{2n+1}$, for more details see~\cite{EkholmEtnyreSullivan2}.

Let $\mathcal C$ be the set of Reeb chords of $\Lambda$. Since $\Lambda$ is generic, $\mathcal C$ is a finite set.
Let $A_{\Lambda}$ be the vector space over $\Z_{2}$ generated by the elements of $\mathcal C$ and let
$\mathcal A_{\Lambda}$ be the unital tensor algebra over $A_{\Lambda}$, i.e.,
\begin{align*}
\mathcal A_{\Lambda}=\bigotimes\limits_{k=0}^{\infty}A_{\Lambda}^{\otimes k}.
\end{align*}
$\mathcal A_{\Lambda}$ is a differential graded algebra whose grading is denoted by $|\cdot|$ and differential is denoted by $\partial_{\Lambda}$. $\mathcal A_{\Lambda}$ is called a Legendrian contact homology differential graded algebra of $\Lambda$. For the definitions of $|\cdot|$ and $\partial_{\Lambda}$ we refer to Section 2 in~\cite{EkholmEtnyreSullivan}.

Note that it is difficult to
use Legendrian contact homology in practical applications, as it is the homology of an
infinite dimensional noncommutative algebra with a nonlinear differential.
One of the ways to extract useful information
from the Legendrian contact homology differential graded algebra is to follow the Chekanov's method of linearization, which
uses an augmentation $\ve:\mathcal A_{\Lambda}\to \Z_{2}$ to produce a finite-dimensional chain complex $LC^{\ve}(\Lambda)$ whose
homology is denoted by $LCH^{\ve}(\Lambda)$~\cite{Chekanov}. More precisely,
$\ve$ is a graded algebra
map $\ve: \mathcal A_{\Lambda}\to \Z_{2}$ that satisfy the following two conditions:
\begin{itemize}
\item[$(1)$] $\ve(1) = 1$;
\item[$(2)$] $\ve\circ \partial_{\Lambda} = 0$.
\end{itemize}
Consider the graded isomorphism $\varphi^{\ve}:\mathcal A_{\Lambda}\to \mathcal A_{\Lambda}$ defined by $\varphi^{\ve}(c)=c+\ve(c)$.
This map defines a new differential $\partial^{\ve}(c):=\varphi^{\ve}\circ\partial_{\Lambda}\circ(\varphi^{\ve})^{-1}(c)$ and $LC^{\ve}(\Lambda):=(A_{\Lambda}, \partial^{\ve}_{1})$, where $\partial^{\ve}_{1}: A_{\Lambda}\to A_{\Lambda}$ is a $1$-component of $\partial^{\ve}$.
We let $LCH_{\ve}(\Lambda)$ be the homology of the dual complex $LC_{\ve}(\Lambda):=Hom(LC^{\ve}(\Lambda),\Z_{2})$.

Following Ekholm~\cite{Ekholmrsftz2}, we observe that exact Lagrangian cobordism between
two Legendrian submanifolds can be used to define a map between the Legendrian contact homology algebras.

In this paper, we establish the following two long exact sequences:

\begin{theorem}\label{homcobcontcohends}
Let $\Lambda_{-}$ and $\Lambda_{+}$ be two closed, orientable Legendrian submanifolds of $\R^{2n+1}$ such that $\emptyset\prec_{L_{\Lambda_{-}}}^{ex}\Lambda_{-}$. Then from the condition $\Lambda_{-}\prec_{L}^{ex}\Lambda_{+}$ it follows that there is the following exact sequence
\begin{align}\label{flexactseqinth}
\rightarrow H_{i}(\Lambda_{-}) \rightarrow H_{i}(L)\oplus LCH_{\ve_{-}}^{n-i+2}(\Lambda_{-})\rightarrow LCH_{\ve_{+}}^{n-i+2}(\Lambda_{+})\rightarrow H_{i-1}(\Lambda_{-}) \rightarrow.
\end{align}
In addition, $\Lambda_{-}\prec_{L}^{ex}\Lambda_{+}$ implies that there is the following exact sequence
\begin{align}\label{slexactseqinth}
\rightarrow LCH_{\ve_{-}}^{n-i+2}(\Lambda_{-})\rightarrow LCH_{\ve_{+}}^{n-i+2}(\Lambda_{+})\rightarrow H_{i} (L, \Lambda_{-})\rightarrow LCH_{\ve_{-}}^{n-i+3}(\Lambda_{-})\rightarrow.
\end{align}
Here $LCH_{\ve_{\pm}}^{i}(\Lambda_{\pm})$ is the linearized Legendrian contact cohomology of $\Lambda_{\pm}$ over $\Z_{2}$, linearized with respect to the augmentation $\ve_{\pm}$.  $\ve_{-}$ is the augmentation induced by $L_{\Lambda_{-}}$, and $\ve_{+}$ is the augmentation induced by $L$ and $\ve_{-}$.
\end{theorem}
We thank Joshua Sabloff and Lisa Traynor for pointing out the way to get the second long exact sequence in Theorem~\ref{homcobcontcohends}.

\subsection{The Thurston-Bennequin invariant}
The Thurston--Bennequin invariant (number) of a closed, orientable, connected Legendrian submanifold
$\Lambda$ of $\R^{2n+1}$ was originally defined by Bennequin~\cite{Bennequin} and independently by Thurston when
$n = 1$, and was generalized to the case when $n\geq 1$ by Tabachnikov~\cite{Tabachnikov}.

Pick an orientation on $\Lambda\subset \R^{2n+1}$. Push $\Lambda$ slightly off of itself
along $R_{\alpha}=\partial_z$ to get another oriented submanifold $\Lambda'$ disjoint from $\Lambda$. The Thurston-Bennequin
invariant of $\Lambda$ is the linking number
\begin{align*}
tb(\Lambda) = lk(\Lambda, \Lambda').
\end{align*}
Note that $tb(\Lambda)$ is independent of the choice of orientation on $\Lambda$
since changing it changes also the orientation of $\Lambda'$.

Our goal is to prove the following theorem:
\begin{theorem}\label{mainthtb}
Let $\Lambda_{-}$ and $\Lambda_{+}$ be two closed, orientable Legendrian submanifolds of $\R^{2n+1}$.
\begin{itemize}
\item[$(1)$] If $n$ is even and $\Lambda_{-}\prec_{L}\Lambda_{+}$, then
\begin{align*}
tb(\Lambda_{+})+tb(\Lambda_{-})=(-1)^{\frac{n}{2}+1}\chi(L).
\end{align*}
\item[$(2)$] If $n$ is odd, $\emptyset\prec_{L_{\Lambda_{-}}}^{ex}\Lambda_{-}$ and
$\Lambda_{-}\prec^{ex}_{L}\Lambda_{+}$, then
\begin{align*}
tb(\Lambda_{+})-tb(\Lambda_{-})=(-1)^{\frac{(n-2)(n-1)}{2}+1}\chi(L).
\end{align*}
\end{itemize}
\end{theorem}

\subsection{Constructions and examples}
In~\cite{Chantraine}, Chantraine described the way to construct Lagrangian cobordisms from Legendrian isotopies of Legendrian knots. We show that the construction of Chantraine works in high dimensions. More precisely, we prove the following:
\begin{proposition}\label{legisotopycobordism}
Let $\Lambda_{-},\Lambda_{+}$ be two closed, orientable Legendrian submanifolds of $\R^{2n+1}$ that are Legendrian isotopic, then
there exists an exact Lagrangian cobordism $L$ such that $\Lambda_{-}\prec^{ex}_{L}\Lambda_{+}$.
\end{proposition}

Front spinning is a procedure to construct a closed, orientable Legendrian submanifold $\Sigma \Lambda\subset\R^{2n+3}$ from a closed, orientable Legendrian submanifold $\Lambda\subset\R^{2n+1}$. It was invented by Ekholm, Etnyre and Sullivan in~\cite{EkholmEtnyreSullivan}. The detailed description of this procedure will be provided in Section~\ref{examples}. We prove the following property of it:
\begin{proposition}\label{frspuncon}
Let $\Lambda_{-},\Lambda_{+}$ be two closed, orientable Legendrian submanifolds of $\R^{2n+1}$. If $\Lambda_{-}\prec_{L}^{lag}\Lambda_{+}$, then there exists a Lagrangian cobordism $\Sigma L$ such that
$\Sigma \Lambda_{-}\prec_{\Sigma L}^{lag}\Sigma \Lambda_{+}$. In addition, if $\Lambda_{-}\prec_{L}^{ex}\Lambda_{+}$, then there exists an exact Lagrangian cobordism $\Sigma L$ such that
$\Sigma \Lambda_{-}\prec_{\Sigma L}^{ex}\Sigma \Lambda_{+}$.
\end{proposition}

Finally, we apply Proposition~\ref{frspuncon} to the exact Lagrangian cobordisms from~\cite{EkholmHondaKalman} and construct exact Lagrangian cobordisms
between the non-isotopic Legendrian tori described in~\cite{EkholmEtnyreSullivan}.
\begin{proposition}\label{pairstoricobtnotisot}
There are infinitely many pairs of exact Lagrangian cobordant and not pairwise
Legendrian isotopic Legendrian $n$-tori in $\R^{2n+1}$.
\end{proposition}

\section{Proof of Theorem~\ref{homcobcontcohends}}
In this section, we prove the existence of the long exact sequences described in Theorem~\ref{homcobcontcohends}.
We first construct an exact Lagrangian filling of $\Lambda_{+}$.

Since $\Lambda_{-}$ is connected, and $L$, $L_{\Lambda_{-}}$ are exact Lagrangian cobordisms in the symplectization of $\R^{2n+1}$
such that ($-\infty$)-boundary of $L$, which is $\Lambda_{-}$, agrees with
($+\infty$)-boundary of $L_{\Lambda_{-}}$, then $L$ and $L_{\Lambda_{-}}$ can
be joined to the exact Lagrangian cobordism $L_{\Lambda_{+}}$ in the symplectization of $\R^{2n+1}$, where $L_{\Lambda_{+}}$ is obtained
by gluing the positive end of $L_{\Lambda_{-}}$ to the negative
end of $L$. Since the $-\infty$-boundary of $L_{\Lambda_{-}}$ is empty, the $-\infty$-boundary of $L_{\Lambda_{+}}$ is also empty.

We now use the Mayer-Vietoris long exact sequence for $L_{\Lambda_{-}}, L\subset L_{\Lambda_{+}}$. We possibly extend $L_{\Lambda_{-}}$ and $L$ in such a way that $L_{\Lambda_{-}}\cap L$ is diffeomorphic to $\R\times \Lambda_{-}$. Hence, the Mayer-Vietoris long exact sequence can be written as
\begin{align*}
\rightarrow H_{i}(\R\times \Lambda_{-})\rightarrow H_{i}(L)\oplus H_{i}(L_{\Lambda_{-}})\rightarrow H_{i}(L_{\Lambda_{+}})\rightarrow H_{i-1}(\R\times \Lambda_{-})\rightarrow.
\end{align*}
Now we note that $H_{i}(\R\times \Lambda_{-})\simeq H_{i}(\Lambda_{-})$ for all $i$. Hence, we can rewrite the Mayer-Vietoris long exact sequence as
\begin{align}\label{mvlexseqsobfills}
\rightarrow H_{i}(\Lambda_{-}) \rightarrow H_{i}(L)\oplus H_{i}(L_{\Lambda_{-}})\rightarrow H_{i}(L_{\Lambda_{+}})\rightarrow H_{i-1}(\Lambda_{-}) \rightarrow.
\end{align}

We now remind the reader of the following fact described by Ekholm in~\cite{Ekholm}, which comes from certain observations of Seidel in wrapped Floer homology~\cite{AbouzaidSeidel}, \cite{FukayaSeidelSmith}:
\begin{fact}[\cite{Ekholm}]\label{legconthomhomfil}
Let $\Lambda$ be a closed, orientable, connected, chord generic Legendrian submanifold of $\R^{2n+1}$ and $\emptyset\prec^{ex}_{L_{\Lambda}}\Lambda$. Then
\begin{align}\label{agagainsignlchcohregfil}
H_{n-i+2}(L_{\Lambda})\simeq LCH_{\ve}^{i}(\Lambda).
\end{align}
Here $\ve$ is the augmentation induced by $L_{\Lambda}$.
\end{fact}
For the definition of the augmentation induced by a filling we refer to Section 3 in~\cite{Ekholmrsftz2}.
Observe that Ekholm in~\cite{Ekholm} provided a fairly complete sketch of proof of Fact~\ref{legconthomhomfil}.

We change the indices in Formula~\ref{agagainsignlchcohregfil} and write it as
\begin{align}\label{agaaagainsignlchcohregfil}
H_{i}(L_{\Lambda_{\pm}})\simeq LCH_{\ve_{\pm}}^{n-i+2}(\Lambda_{\pm}).
\end{align}
Using Formula~\ref{agaaagainsignlchcohregfil}, we rewrite Mayer-Vietoris long exact sequence~\ref{mvlexseqsobfills} as
\begin{align}\label{finversionleslchnl}
\rightarrow H_{i}(\Lambda_{-}) \rightarrow H_{i}(L)\oplus LCH_{\ve_{-}}^{n-i+2}(\Lambda_{-})\rightarrow LCH_{\ve_{+}}^{n-i+2}(\Lambda_{+})\rightarrow H_{i-1}(\Lambda_{-}) \rightarrow.
\end{align}

We now write the long exact sequence for the pair $(L_{\Lambda_{-}}, L_{\Lambda_{+}})$
\begin{align}\label{leqpairfst}
\rightarrow H_{i}(L_{\Lambda_{-}})\rightarrow H_{i}(L_{\Lambda_{+}})\rightarrow H_{i} (L_{\Lambda_{+}}, L_{\Lambda_{-}})\rightarrow H_{i-1}(L_{\Lambda_{-}})\rightarrow.
\end{align}
Using Formula~\ref{agaaagainsignlchcohregfil} and the excision theorem for $L_{\Lambda_{+}}, L\subset  L_{\Lambda_{+}}$, we
write long exact sequence~\ref{leqpairfst} as
\begin{align}\label{seclchgoodone}
\rightarrow LCH_{\ve_{-}}^{n-i+2}(\Lambda_{-})\rightarrow LCH_{\ve_{+}}^{n-i+2}(\Lambda_{+})\rightarrow H_{i} (L, \Lambda_{-})\rightarrow LCH_{\ve_{-}}^{n-i+3}(\Lambda_{-})\rightarrow.
\end{align}
This  finishes the proof of Theorem~\ref{homcobcontcohends}.

\begin{remark}\label{factorlinconthom}
Note that under the conditions of Theorem~\ref{homcobcontcohends}, if $H_{i}(\Lambda_{-})=H_{i-1}(\Lambda_{-})=0$ for some $i$, say when $\Lambda_{-}=S^{n}$ and $i,i-1\neq 0,n$, then long exact sequence~\ref{finversionleslchnl} implies that
\begin{align*}
LCH_{\ve_{+}}^{n-i+2}(\Lambda_{+})\simeq H_{i}(L)\oplus LCH_{\ve_{-}}^{n-i+2}(\Lambda_{-}).
\end{align*}
Hence, for such $i$ we get that
\begin{align*}
H_{i}(L)\simeq LCH_{\ve_{+}}^{n-i+2}(\Lambda_{+})/LCH_{\ve_{-}}^{n-i+2}(\Lambda_{-}).
\end{align*}
\end{remark}

\begin{remark}
Note that we can rewrite long exact sequences~\ref{finversionleslchnl} and~\ref{seclchgoodone} using the relative symplectic field theory of $((\R\times \R^{2n+1}, d(e^t\alpha)), L_{\Lambda_{\pm}})$, since
\begin{align}\label{sftlch}
E_{1}^{i}((\R\times \R^{2n+1}, d(e^t\alpha)), L_{\Lambda_{\pm}})\simeq LCH_{\ve_{\pm}}^{i}(\Lambda_{\pm})
\end{align}
over $\Z_{2}$.
For the definition of the relative symplectic field theory we refer to~\cite{Ekholmrsftz2}, for the details about the isomorphism described in Formula~\ref{sftlch} we refer to~\cite{Ekholm} (we observe that since $L_{\Lambda_{\pm}}$ are
connected, the associated spectral sequences have only one level).
\end{remark}

\section{Proof of Theorem~\ref{mainthtb}}\label{tbtheorfordefconthom}
Let $n$ be even.
We first recall the following proposition from~\cite{Eliashberg}:
\begin{proposition}[\cite{Eliashberg}]\label{propreprestbaseulchar}
Let $\Lambda$ be a closed, orientable, connected, chord generic Legendrian submanifold  of $\R^{2n+1}$, where $n$ is even. Then
\begin{align*}
tb(\Lambda)=(-1)^{\frac{n}{2}+1}\frac{1}{2}\chi(\Lambda).
\end{align*}
\end{proposition}

We now note that
\begin{align}\label{eulcharboundint}
\chi(\partial L)=2\chi(L).
\end{align}
Equation~\ref{eulcharboundint} holds because the
Euler characteristic of an even-dimensional boundary is twice the Euler characteristic of its bounded
manifold, see Chapter 21 in~\cite{May}. We now observe that $\partial L = \Lambda_{+} \sqcup \Lambda_{-}$ and hence from Equation~\ref{eulcharboundint} we get that
\begin{align}\label{eulcharboundintalmost}
2\chi(L)=\chi(\partial L)=\chi(\Lambda_{+})+\chi(\Lambda_{-}).
\end{align}
Then we use Proposition~\ref{propreprestbaseulchar} and rewrite Equation~\ref{eulcharboundintalmost} as
\begin{align}\label{eulcharboundintalmostalmost}
2\chi(L)=\chi(\Lambda_{+})+\chi(\Lambda_{-})=2(-1)^{-\frac{n}{2}-1}(tb(\Lambda_{+})+tb(\Lambda_{-})).
\end{align}
From Equation~\ref{eulcharboundintalmostalmost} it follows that
\begin{align}
tb(\Lambda_{+})+tb(\Lambda_{-})=(-1)^{\frac{n}{2}+1}\chi(L).
\end{align}
This finishes the proof of Theorem~\ref{mainthtb} in the case when $n$ is even.

We now prove $(2)$.
First we provide another definition of the Thurston-Bennequin number. The author learned about the alternative definition from~\cite{EkholmEtnyreSullivan2}.

Let $\Lambda$ be a closed, orientable, connected, chord generic Legendrian submanifold of $\R^{2n+1}$ and let $c$ be a Reeb chord of $\Lambda$ with end points $a$ and $b$ such that $z(a) > z(b)$.
We define $V_{a}:= d\Pi(T_{a}\Lambda)$ and
$V_{b}:= d\Pi(T_{b}\Lambda)$. Given an orientation on $\Lambda$, $V_{a}$ and $V_{b}$ are oriented $n$-dimensional transverse
subspaces of $\R^{2n}$. If the orientation of $V_{a}\oplus V_{b}$ agrees with that of $\R^{2n}$, then we say that the sign of $c$, we denote it by
$sign(c)$,  is $+1$, otherwise we say that it is $-1$. Then
\begin{align}\label{onedefntb}
tb(\Lambda) = \sum\limits_{c} sign(c),
\end{align}
where the sum is taken over all Reeb chords $c$ of $\Lambda$.

Using Formula~\ref{onedefntb}, the following proposition was proven in~\cite{EkholmEtnyreSullivan}:
\begin{proposition}[\cite{EkholmEtnyreSullivan}]\label{tbeulcharlincont}
If $\Lambda\subset \R^{2n+1}$ is a closed, orientable, connected, chord generic Legendrian submanifold, then
\begin{align*}
tb(\Lambda)=(-1)^{\frac{(n-2)(n-1)}{2}}\sum\limits_{c\in \mathcal C} (-1)^{|c|}.
\end{align*}
\end{proposition}

We now construct an exact Lagrangian filling of $\Lambda_{+}$. We do it the same way as in the proof of Theorem~\ref{homcobcontcohends}, namely
$L_{\Lambda_{+}}$ is obtained by gluing the positive end of $L_{\Lambda_{-}}$ to the negative
end of $L$ in the symplectization of $\R^{2n+1}$.

Taking Euler characteristics of the long exact sequence~\ref{slexactseqinth} and using Proposition~\ref{tbeulcharlincont}, we get

\begin{align}\label{finaltbodd}
tb (\Lambda_{+}) - tb (\Lambda_{-})=(-1)^{\frac{(n-2)(n-1)}{2}+1}\chi(L).
\end{align}
This finishes the proof of Theorem~\ref{mainthtb} when $n$ is odd.
\begin{remark}
Note that when $n=1$ Equation~\ref{finaltbodd} can be written  as
\begin{align*}
tb (\Lambda_{+}) - tb (\Lambda_{-})=-\chi(L),
\end{align*}
which coincides with the formula from Theorem 1.2 in~\cite{Chantraine}.
\end{remark}

\begin{remark}
Observe that the condition of Theorem~\ref{mainthtb} in the case when $n$ is odd is much stronger than the condition of Theorem~\ref{mainthtb} in the case when $n$ is even.
If $n$ is even, $\emptyset\prec_{L_{\Lambda_{-}}}^{ex}\Lambda_{-}$ and
$\Lambda_{-}\prec^{ex}_{L}\Lambda_{+}$, then taking Euler characteristics of the long exact sequence~\ref{slexactseqinth} and using Proposition~\ref{tbeulcharlincont} we get that
\begin{align*}
tb(\Lambda_{+})+tb(\Lambda_{-})=(-1)^{\frac{n}{2}+1}\chi(L).
\end{align*}
\end{remark}

Note that the proof of Theorem~\ref{mainthtb} can be easily modified to become a proof of the following remark:
\begin{remark}\label{thmfiltb}
Let $\Lambda$ be a closed, orientable Legendrian submanifold of $\R^{2n+1}$.
\begin{itemize}
\item[$(1)$] If $n$ is even and $\emptyset\prec_{L_{\Lambda}}\Lambda$, then
\begin{align*}
tb(\Lambda)=(-1)^{\frac{n}{2}+1}\chi(L_{\Lambda}).
\end{align*}
\item[$(2)$] If $n$ is odd and $\emptyset\prec_{L_{\Lambda}}^{ex}\Lambda$, then
\begin{align*}
tb(\Lambda)=(-1)^{\frac{(n-2)(n-1)}{2}+1}\chi(L_{\Lambda}).
\end{align*}
\end{itemize}
\end{remark}

\section{Examples}\label{examples}
In this section, we describe a few examples of Lagrangian cobordisms. These examples are based on the works of Chantraine~\cite{Chantraine}, Ekholm, Etnyre and Sullivan~\cite{EkholmEtnyreSullivan}, and Ekholm, Honda and K\'{a}lm\'{a}n~\cite{EkholmHondaKalman}. For the constructions of Lagrangian cobordisms based
on the generating families technique we refer to Sabloff and Traynor~\cite{SabloffTraynor2}.

\example[proof of Proposition~\ref{legisotopycobordism}]
Let $\Lambda_{-}, \Lambda_{+}\subset \R^{2n+1}$ be two closed, orientable Legendrian submanifolds which are Legendrian isotopic.
This means that there
is a smooth isotopy of a closed manifold $\Lambda$ to $\R^{2n+1}$ given by $\varphi: \Lambda\times [0,1]\to \R^{2n+1}$ such that $\Lambda_{\nu}:=\varphi(\Lambda,{\nu})$ is Legendrian for all $\nu\in[0,1]$, $\Lambda_{-}=\Lambda_{0}$ and $\Lambda_{+}=\Lambda_{1}$. We now construct $L$ such that $\Lambda_{-}\prec^{ex}_{L} \Lambda_{+}$. Observe that in the construction below one can omit the assumption that $\Lambda_{-}, \Lambda_{+}, L$ are connected.
In the case of Legendrian knots in $\R^{3}$, the construction of $L$ was described by Chantraine, see Theorem~1.1 in~\cite{Chantraine}.
In our case, the construction of Chantraine can be described in the following way:
\begin{itemize}
\item[(1)] We note that $\R\times \Lambda_{-}$ is an exact Lagrangian submanifold of $(\R\times \R^{2n+1}, d(e^t\alpha))$.
\item[(2)]  Theorem~2.6.2 from~\cite{Geiges} implies that there is a compactly supported one-parameter family
of contactomorphisms $f_{\nu}$ which realizes the isotopy $(\Lambda_{\nu})_{\nu\in[0,1]}$.
\item[(3)] Proposition~2.2 from~\cite{Chantraine} implies that a contactomorphism of $\R^{2n+1}$ lifts to a Hamiltonian diffeomorphism of the symplectization
$(\R\times \R^{2n+1}, d(e^t\alpha))$.
\item[(4)] Let $H$ be a Hamiltonian on
$\R\times \R^{2n+1}$ whose flow realizes the lifts of $f_{\nu}$'s. The existence of $H$ follows from (3). Following Chantraine, we construct
\begin{align*}
H':\R\times \R^{2n+1}\times[0,1]\to\R
\end{align*}
such that
\begin{align*}
H'(t,x,\nu)=\left \{ \begin{array}{ll} H(t,x,\nu), & \mbox{for}\ t>T; \\
0, & \mbox{for}\ t<-T.
\end{array} \right.
\end{align*}
Here $T\gg 0$.
\item[(5)] Let $\phi^{\nu}$ be the Hamiltonian flow of $H'$. We now observe that $\phi^{1}(\R\times \Lambda_{-})$ coincides with $\R\times \Lambda_{-}$ near $-\infty$ and with $\R\times \Lambda_{+}$ near $\infty$.
\item[(6)] Since $\R\times \Lambda_{-}$ is exact and $\phi^{1}$ is a Hamiltonian diffeomorphism, $L:=\phi^{1}(\R\times \Lambda_{-})$ is exact.
\end{itemize}
This finishes the proof of Proposition~\ref{legisotopycobordism}.

\begin{remark}
Note that Eliashberg and Gromov in~\cite{EliashbergGromov} provided another proof of the fact that Legendrian isotopy implies Lagrangian cobordism.
\end{remark}

\example[proof of Proposition~\ref{frspuncon}]\label{exfrspinex}
The following construction is based on the front spinning method invented by Ekholm, Etnyre and Sullivan in~\cite{EkholmEtnyreSullivan}.

First we recall the notion of front projection.
{\em Front projection} is a map $\Pi_{F}$ from $\R^{2n+1}$ to $\R^{n+1}$ defined by
\begin{align*}
\Pi_{F}(x_{1},y_{1},\dots,x_{n},y_{n},z)=(x_{1},x_{2},\dots,x_{n},z).
\end{align*}
Let $\Lambda$ be a closed, orientable Legendrian submanifold of $\R^{2n+1}$ parameterized by $f_{\Lambda}:\Lambda\to \R^{2n+1}$ and we write
\begin{align*}
f_{\Lambda}(p)=(x_{1}(p),y_{1}(p),\dots,x_{n}(p),y_{n}(p),z(p))
\end{align*}
for $p\in \Lambda$.
The front projection of $\Lambda$ is parameterized by $\Pi_{F}\circ f_{\Lambda}$ and we have
\begin{align*}
\Pi_{F}\circ f_{\Lambda}(p)=(x_{1}(p),x_{2}(p),\dots,x_{n}(p),z(p)).
\end{align*}
Without loss of generality we can assume that $x_{1}(p)>0$ for all $p\in \Lambda$. We now embed $\R^{n+1}$ to $\R^{n+2}$ via
\begin{align*}
(x_{1},\dots,x_{n},z)\to (x_{0}=0,x_{1},\dots,x_{n},z)
\end{align*}
and construct the suspension
of $\Lambda$, we denote it by $\Sigma \Lambda$, such that $\Pi_{F}(\Sigma \Lambda)$ is obtained from $\Pi_{F}(\Lambda)$ by rotating it around the subspace
$x_{0}=x_{1}=0$. $\Pi_{F}(\Sigma \Lambda)$ can be parameterized by $(x_{1}(p)\sin\theta,x_{1}(p)\cos\theta,x_{2}(p),\dots,x_{n}(p),z(p))$ with $\theta\in S^1$
and is the front projection of a Legendrian embedding $\Lambda\times S^1\to \R^{2n+3}$. For the properties of $\Sigma \Lambda$ we refer to Lemma 4.16 in~\cite{EkholmEtnyreSullivan}.

Let $\Lambda_{-}, \Lambda_{+}$ be two closed, orientable Legendrian submanifolds of $\R^{2n+1}$ such that
\begin{align}\label{poslegtotransl}
\Lambda_{\pm}\subset \{(x_{1},y_{1},\dots,x_{n},y_{n},z)\in \R^{2n+1}\ |\ x_{1}>0\}
\end{align}
and
$\Lambda_{-}\prec^{lag}_{L}\Lambda_{+}$. Let $L$ be parameterized by
$f_{L}:L\to \R^{2n+2}$
\begin{align*}
f_{L}(p)=(t(p),x_{1}(p),y_{1}(p),\dots,x_{n}(p),y_{n}(p),z(p)).
\end{align*}
Without loss of generality we assume that $x_{1}(p)>0$ for all $p$ (Formula~\ref{poslegtotransl} implies that
\begin{align*}
\{f_{L}(p)\ |\ x_{1}(p)\leq 0\}
\end{align*}
is compact and we can translate $L$ so that $x_{1}(p)>0$ for all $p$).
Then we construct a Lagrangian cobordism from $\Sigma \Lambda_{-}$ to $\Sigma \Lambda_{+}$ that we call $\Sigma L$. We define $\Sigma L$ to be parameterized by
\begin{align*}
f_{\Sigma L}: L\times S^1\to \R\times \R^{2n+3}
\end{align*}
with
\begin{align*}
f_{\Sigma L}(p,\theta)=(t(p),x_{1}(p)\sin\theta,y_{1}(p)\sin\theta,x_{1}(p)\cos\theta,y_{1}(p)\cos\theta,x_{2}(p),\dots, z(p)).
\end{align*}
Here $p\in L$ and $\theta\in S^1$.

We now show that $\Sigma L$ is really a Lagrangian cobordism from $\Sigma \Lambda_{-}$ to $\Sigma \Lambda_{+}$.
Let
\begin{align*}
&\Lambda_{+}^{T_{L}}:=\{(x_{0},\dots ,y_{n},z)\ |\ (T_{L},x_{0},\dots,y_{n},z)\in f_{\Sigma L}(\Sigma L)\cap (\{T_{L}\}\times\R^{2n+3})\}\ \
\mbox{and}\\
&\Lambda_{-}^{T_{L}}:=\{(x_{0},\dots ,y_{n},z)\ |\ (-T_{L},x_{0},\dots,y_{n},z)\in f_{\Sigma L}(\Sigma L)\cap (\{-T_{L}\}\times \R^{2n+3})\}.
\end{align*}
We now note that from the definition of $T_{L}$ it follows that
\begin{align*}
&f_{\Sigma L}(\Sigma L)\cap ([T_{L},\infty)\times \R^{2n+3})=[T_{L},\infty)\times \Lambda_{+}^{T_{L}}\ \
\mbox{and}\\
&f_{\Sigma L}(\Sigma L)\cap ((-\infty,-T_{L}]\times \R^{2n+3})=(-\infty,-T_{L}]\times \Lambda_{-}^{T_{L}}.
\end{align*}
In addition, we observe that $\Lambda_{\pm}^{T_{L}}\subset \R^{2n+3}$ can be parameterized by
\begin{align*}
f_{\Lambda_{\pm}^{T_{L}}}:\Lambda_{\pm}\times S^1\to \R^{2n+3}
\end{align*}
such that
\begin{align*}
f_{\Lambda_{\pm}^{T_{L}}}(p,\theta)=(x_{1}(p)\sin\theta,y_{1}(p)\sin\theta,x_{1}(p)\cos\theta,y_{1}(p)\cos\theta,x_{2}(p),\dots, z(p)).
\end{align*}
Here $p\in \Lambda_{\pm}\subset \partial L$ and $\theta\in S^1$.
We now prove that $\Lambda_{\pm}^{T_{L}}$ coincides with $\Sigma \Lambda_{\pm}$. It is clear that $\Pi_{F}(\Lambda_{\pm}^{T_{L}})=\Pi_{F}(\Sigma\Lambda_{\pm})$.
It remains to prove that $\Lambda_{\pm}^{T_{L}}$ is a Legendrian submanifold of $\R^{2n+3}$.

It is easy to see that
\begin{align}\label{pullbcontformworks}
&f_{\Lambda_{\pm}^{T_{L}}}^{\ast}(dz-\sum\limits_{i=0}^{n}y_{i}dx_{i}) = dz(p)-\sum\limits_{i=2}^{n}y_{i}(p)dx_{i}(p)\\&-y_{1}(p)(\sin^2 \theta + \cos^2 \theta)dx_{1}(p)
+(y_{1}(p)x_{1}(p)\sin\theta\cos\theta - y_{1}(p)x_{1}(p)\sin\theta\cos\theta)d\theta\nonumber.
\end{align}
Since $\Lambda_{\pm}$ is Legendrian submanifold of $\R^{2n+1}$ and hence $f_{\Lambda_{\pm}}^{\ast}(dz-\sum_{i=1}^{n}y_{i}dx_{i})=0$, we have that
\begin{align}\label{prevnewconform}
y_{1}(p)dx_{1}(p)=dz(p)-\sum\limits_{i=2}^{n}y_{i}(p)dx_{i}(p).
\end{align}
Hence, Formulas~\ref{pullbcontformworks} and~\ref{prevnewconform} imply that
\begin{align}\label{pullbcontformst}
f_{\Lambda_{\pm}^{T_{L}}}^{\ast}(dz-\sum_{i=0}^{n}y_{i}dx_{i})=0.
\end{align}
Since
\begin{align*}
f_{\Lambda_{\pm}}(p):=(x_{1}(p),\dots,y_{n}(p),z(p))
\end{align*}
with $p\in \Lambda_{\pm}\subset \partial L$ is a parametrization of an embedded submanifold of dimension $n$, and $x_{1}(p)>0$ for $p\in \Lambda_{\pm}\subset \partial L$, one easily sees that
\begin{align*}
f_{\Lambda_{\pm}^{T_{L}}}(p)=(x_{1}(p)\sin\theta,y_{1}(p)\sin\theta,x_{1}(p)\cos\theta,y_{1}(p)\cos\theta,x_{2}(p),\dots, z(p))
\end{align*}
with $p\in \Lambda_{\pm}$, $\theta\in S^1$ is a parametrization of an embedded submanifold of dimension $n+1$.
Thus, using Formula~\ref{pullbcontformst} we see that $\Lambda_{\pm}^{T_{L}}$ is an embedded Legendrian submanifold of $\R^{2n+3}$ whose front projection coincides with $\Pi_{F}(\Sigma \Lambda_{\pm})$. Thus, we get that
$\Lambda_{\pm}^{T_{L}}=\Sigma \Lambda_{\pm}$.

We now note that
\begin{align}\label{longpullbsymplnzero}
&f_{\Sigma L}^{\ast}(d(e^t(dz-\sum\limits_{i=0}^{n}y_{i}dx_{i}))) = e^t(dt(p)\wedge dz(p)-\sum\limits_{i=2}^{n}dy_{i}(p)\wedge dx_{i}(p)\\ &-\sum\limits_{i=2}^{n}y_{i}(p)dt(p)\wedge dx_{i}(p)-(y_{1}(p)(\sin^{2}\theta + \cos^{2}\theta) dt(p)\wedge dx_{1}(p)\nonumber\\
&+(\sin^{2}\theta+\cos^{2}\theta)dy_{1}(p)\wedge dx_{1}(p)+(\sin^{2}\theta+\cos^{2}\theta)x_{1}(p)y_{1}(p)d\theta\wedge d\theta\nonumber\\
&+(y_{1}(p)x_{1}(p)\sin \theta \cos \theta-y_{1}(p)x_{1}(p)\sin \theta \cos \theta ) dt(p)\wedge d\theta\nonumber\\
&+(y_{1}(p)\sin\theta\cos\theta-y_{1}(p)\sin\theta\cos\theta)d\theta\wedge dx_{1}(p)\nonumber\\
&+(x_{1}(p)\sin\theta\cos\theta-x_{1}(p)\sin\theta\cos\theta)dy_{1}(p)\wedge d\theta ))\nonumber.
\end{align}
In addition, observe that
\begin{align}\label{symplecnfrompreviouscoordtonew}
&e^t(dt(p)\wedge dz(p)-\sum\limits_{i=2}^{n}dy_{i}(p)\wedge dx_{i}(p) -\sum\limits_{i=2}^{n}y_{i}(p)dt(p)\wedge dx_{i}(p))\\&=e^t(y_{1}(p)dt(p)\wedge dx_{1}(p)+dy_{1}(p)\wedge dx_{1}(p))\nonumber.
\end{align}
Hence, Formulas~\ref{longpullbsymplnzero} and~\ref{symplecnfrompreviouscoordtonew} imply that
\begin{align}\label{pullbsympln}
f_{\Sigma L}^{\ast}(d(e^t(dz-\sum^{n}_{i=0}y_{i}dx_{i})))=0.
\end{align}
Since
\begin{align*}
f_{L}(p)=(t(p),x_{1}(p),y_{1}(p),\dots,x_{n}(p),y_{n}(p),z(p)),
\end{align*}
where $p\in L$, is a parametrization of an embedded cobordism of dimension $n+1$ and $x_{1}(p)>0$ for $p\in L$, one easily sees that
\begin{align*}
f_{\Sigma L}(p,\theta)=(t(p),x_{1}(p)\sin\theta,y_{1}(p)\sin\theta,x_{1}(p)\cos\theta,y_{1}(p)\cos\theta,x_{2}(p),\dots, z(p)),
\end{align*}
where $p\in L$ and $\theta\in S^1$, is a parametrization of an embedded cobordism of dimension $n+2$. Hence, we use Formula~\ref{pullbsympln} and see that $\Sigma L$ is really an embedded Lagrangian cobordism from $\Sigma \Lambda_{-}$ to $\Sigma \Lambda_{+}$.

We now assume that $\Lambda_{-}\prec^{ex}_{L}\Lambda_{+}$. Then there is a function $h_{L}\in C^{\infty}(f_{L}(L),\R)$ such that
\begin{align*}
dh_{L}=e^t(dz-\sum^{n}_{i=1}y_{i}dx_{i}).
\end{align*}
From the calculation similar to  Formula~\ref{pullbcontformworks} it follows that
\begin{align}\label{exformcontsympl}
f_{\Sigma L}^{\ast}(e^t(dz-\sum^{n}_{i=0}y_{i}dx_{i}))=e^{t(p)}(dz(p)-\sum^{n}_{i=1}y_{i}(p)dx_{i}(p)).
\end{align}
Since $f_{\Sigma L}$ is an embedding, we can define $h_{\Sigma L}\in C^{\infty}(f_{\Sigma L}(\Sigma L),\R)$ by setting
\begin{align*}
(f_{\Sigma L}^{\ast}h_{\Sigma L})(p,\theta):=(f_{L}^{\ast}h)(p).
\end{align*}
Hence, we use Formula~\ref{exformcontsympl} and get that
\begin{align}\label{exactnesslongformpushed}
d(f_{\Sigma L}^{\ast}h_{\Sigma L})=e^{t(p)}(dz(p)-\sum^{n}_{i=1}y_{i}(p)dx_{i}(p)))=f_{\Sigma L}^{\ast}(e^t(dz-\sum^{n}_{i=0}y_{i}dx_{i})).
\end{align}
Therefore, since $f_{\Sigma L}$ is an embedding, Formula~\ref{exactnesslongformpushed} implies that
\begin{align*}
d(h_{\Sigma L})=e^t(dz-\sum^{n}_{i=0}y_{i}dx_{i}).
\end{align*}
Hence, $\Sigma L$ is an exact Lagrangian cobordism.
This finishes the proof of Proposition~\ref{frspuncon}.

Note that the proof of Proposition~\ref{frspuncon} can be easily modified to become a proof of the following remark:
\begin{remark}\label{thmfiltb}
Let $\Lambda$ be a closed, orientable Legendrian submanifolds of $\R^{2n+1}$. If $\emptyset\prec_{L_{\Lambda}}^{lag}\Lambda$, then there exists a Lagrangian filling $L_{\Sigma\Lambda}$ such that
$\emptyset \prec_{L_{\Sigma\Lambda}}^{lag}\Sigma \Lambda$. In addition, if $\emptyset\prec_{L_{\Lambda}}^{ex}\Lambda$, then there exists an exact Lagrangian filling $L_{\Sigma\Lambda}$ such that
$\emptyset\prec_{L_{\Sigma\Lambda}}^{ex}\Sigma \Lambda$.
\end{remark}

Before we discuss the next example, we briefly recall a few facts about exact Lagrangian cobordisms between Legendrian knots in $\R^{3}$.

\begin{figure}[h]
\begin{center}
\includegraphics[width=250pt]{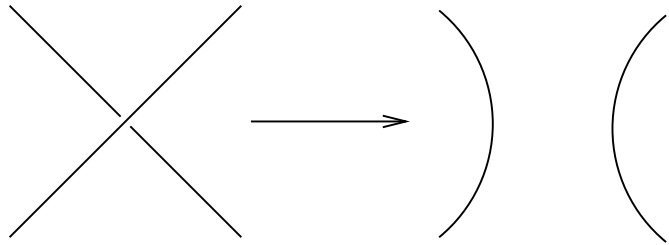}
\caption{The $0$-resolution on the Lagrangian projection.}\label{0resnlagproj}
\end{center}
\end{figure}

\begin{theorem}[\cite{EkholmHondaKalman,EkholmHondaKalman2}]\label{exLagrcobrealization}
There exists an exact Lagrangian cobordism for the following:
\begin{itemize}
\item[(1)] Legendrian isotopy,
\item[(2)] $0$-resolution at a contractible crossing in the Lagrangian projection,
\item[(3)] capping off a $tb=-1$ unknot with a disk.
\end{itemize}
\end{theorem}
See Figure~\ref{0resnlagproj} for the $0$-resolution on the Lagrangian projection.

Following Ekholm, Honda and K\'{a}lm\'{a}n, we
say that a \emph{contractible crossing} of $\Lambda$ is a crossing so that $z_{1}-z_{0}$ can be shrunk to
zero without affecting the other crossings. (Here $z_{1}$ is the $z$-coordinate on
the upper strand and $z_{0}$ is the $z$-coordinate on the lower strand.)
\begin{remark}
Note that Chantraine in~\cite{Chantraine} proved the first part of Theorem~\ref{exLagrcobrealization}.
\end{remark}

\begin{remark}
Note that the second part of Theorem~\ref{exLagrcobrealization} can be proven using the model from Section 3.3 in~\cite{Rizell}.
\end{remark}

\begin{conjecture}[\cite{EkholmHondaKalman,EkholmHondaKalman2}]
If $\emptyset\prec^{ex}_{L_{\Lambda}}\Lambda$, then $L_{\Lambda}$ is obtained
by stacking exact Lagrangians cobordisms described in
Theorem~\ref{exLagrcobrealization}.
\end{conjecture}

\begin{figure}[h]
\begin{center}
\includegraphics[width=250pt]{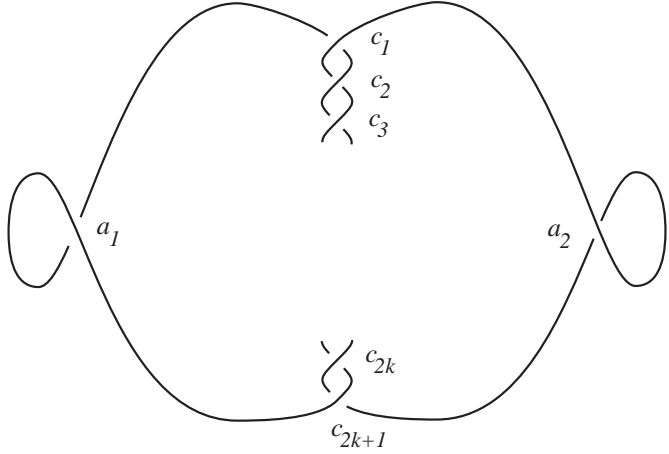}
\caption{The knot $T_{2k+1}$, cf Figure~13 in~\cite{EkholmEtnyreSullivan}.}\label{nonlegisknotsexamples}
\end{center}
\end{figure}

\example[Proof of Proposition~\ref{pairstoricobtnotisot}]
We now use Example~\ref{exfrspinex} to get infinitely many pairs of exact Lagrangian cobordant and not pairwise
Legendrian isotopic Legendrian $n$-tori in $\R^{2n+1}$.
We first recall that Theorem~\ref{exLagrcobrealization} says that $0$-resolution at a contractible crossing in the Lagrangian projection can be realized as an exact Lagrangian cobordism.
Let $T_{2k+1}$ be the Legendrian torus knot from Example~4.18 in~\cite{EkholmEtnyreSullivan}, see Figure~\ref{nonlegisknotsexamples} for the Lagrangian projection of $T_{2k+1}$. One observes that all the crossings in the middle part of the Lagrangian projection are
contractible, see \cite{EkholmHondaKalman2} for the case of $T_{3}$, and hence one can get $T_{2k-1}$ from $T_{2k+1}$ by contracting $c_{2k+1}$ and then $c_{2k}$.  Let $L^{2k+1}_{2k}$ be an exact Lagrangian cobordism which corresponds to the $0$-resolution at $c_{2k+1}$ and let $L^{2k}_{2k-1}$ be an exact Lagrangian cobordism from $T_{2k-1}$ to $T_{2k}$ which corresponds to the resolution of $c_{2k}$. Then we stack $L^{2k+1}_{2k}$ and $L^{2k}_{2k-1}$ and get an exact Lagrangian cobordism that we call $L_{2k-1}^{2k+1}$ such that $T_{2k-1}\prec_{L^{2k+1}_{2k-1}}^{ex}T_{2k+1}$. If we stack $L_{2i-1}^{2i+1}$'s we get an exact Lagrangian cobordism $L^{2k+1}_{2j+1}$ such that $T_{2j+1}\prec_{L^{2k+1}_{2j+1}}^{ex}T_{2k+1}$ for $k>j$.
We use the construction described in Example~\ref{exfrspinex} and get $\Sigma^{n} T_{2j+1}\prec_{\Sigma^{n} L^{2k+1}_{2j+1}}^{ex}\Sigma^{n}T_{2k+1}$ for $k>j$. We now recall that Ekholm, Etnyre and Sullivan proved that $\Sigma^{n} T_{2j+1}$ is not Legendrian isotopic to $\Sigma^{n} T_{2k+1}$ for $k>j+1$ and $j\in \mathbb N$, see Theorem 4.19 in~\cite{EkholmEtnyreSullivan}.

Hence, we get infinitely many pairs of exact Lagrangian cobordant and not pairwise
Legendrian isotopic Legendrian $n$-tori in $\R^{2n+1}$. This finishes the proof of Proposition~\ref{pairstoricobtnotisot}.

\begin{remark}
Given $n\geq 1$. We observe that Theorem~4.19 from~\cite{EkholmEtnyreSullivan} implies that all the Legendrian $n$-tori from Proposition~\ref{pairstoricobtnotisot} are not distinguished by the classical invariants.
\end{remark}

\section*{Acknowledgements}
The author is deeply grateful to Baptiste Chantraine, Vincent Colin, Olivier Collin, Octav Cornea, Tobias Ekholm, Yakov Eliashberg, John Etnyre, Paolo Ghiggini, Ko Honda, Cl\'{e}ment Hyvrier, Georgios Dimitroglou Rizell, Joshua Sabloff, Lisa Traynor and Vera V\'{e}rtesi for helpful conversations and interest in his work. In addition, the author is grateful to the referee of an earlier version of this
paper for many valuable comments and suggestions.

\end{document}